\documentclass[a4paper,12pt,leqno]{article}
\usepackage{amssymb,latexsym}
\usepackage{mltex}
\usepackage{amsmath}
\usepackage{epsfig}
\usepackage{graphics}
\usepackage{pb-diagram}
 \usepackage{lscape} 
\usepackage{enumerate}
\usepackage{amsthm}
\newtheorem{thm}{Theorem}     
\newtheorem{cor}{Corollary}     
\newtheorem{lem}{Lemma}[section]
     
\newtheorem{rem}{Remark}

\newtheorem*{thrm}{Theorem}
\newcommand{\nablab}{\overline{\nabla}}

\newcommand{\Ric}{\mathrm{Ric}}

\newcommand{\lgra}{\longrightarrow}
\newcommand{\Hinf}{||H||_{\infty}}

\newcommand{\Hkp}{||H_k||_{2p}}
\newcommand{\HHkp}{||H_k||^2_{2p}}

\newcommand{\Binf}{||B||_{\infty}}

\newcommand{\HHp}{||H||^2_{2p}}

\newcommand{\ddiv}{\mathrm{div}\,}

\newcommand{\iinf}{\mathrm{inf}\,}

\newcommand{\ssup}{\mathrm{Sup}\,}
\newcommand{\sscal}{\mathrm{Scal}\,}

\newcommand{\trace}{\mathrm{tr\,}}

\newcommand{\lto}{\ensuremath{\longrightarrow}}

\newcommand{\Ss}{\mathbb{S}}

\newcommand{\R}{\mathbb{R}}

\newcommand{\base}{\{e_1,\ldots ,e_n\}}  
 
\newcommand{\function}[5]
{\begin{eqnarray*}\begin{array}{r@{}ccl}
 #1\;\colon\;  & #2 &\lto & #3 \\[.05cm]  
  & #4 &\longmapsto  & #5 
\end{array}\end{eqnarray*}
}
\newcommand{\beqt}{\begin{equation}}  \newcommand{\eeqt}{\end{equation}}
\newcommand{\bal}{\begin{align}}      \newcommand{\eal}{\end{align}}
\newcommand{\ba}{\begin{array}}      \newcommand{\ea}{\end{array}}
\newcommand{\bc}{\begin{center}}     \newcommand{\ec}{\end{center}}
\newcommand{\be}{\begin{enumerate}}  \newcommand{\ee}{\end{enumerate}}
\newcommand{\beq}{\begin{eqnarray}}  \newcommand{\eeq}{\end{eqnarray}}
\newcommand{\beQ}{\begin{eqnarray*}} \newcommand{\eeQ}{\end{eqnarray*}}
\newcommand{\bi}{\begin{itemize}}    \newcommand{\ei}{\end{itemize}}
\newcommand{\bt}{\begin{tabular}}    \newcommand{\et}{\end{tabular}}
\newcommand{\finpreuve}{\hfill\square\\}
\title{Pinching of the First Eigenvalue of the Laplacian and almost-Einstein Hypersurfaces of the Euclidean Space}
\author{Julien Roth}
\date{}
\begin{document}
\maketitle
\begin{center}
Institut \'Elie Cartan, UMR 7502\\
Nancy-Universit\'e, CNRS, INRIA\\
B.P. 239, 54506 Vand\oe uvre l\`es Nancy Cedex, France
\end{center}
\begin{center}
roth@iecn.u-nancy.fr
\end{center}
\begin{abstract}
In this paper, we prove new pinching theorems for the first eigenvalue $\lambda_1(M)$ of the Laplacian on compact hypersurfaces of the Euclidean space. These pinching results are associated with the upper bound for $\lambda_1(M)$ in terms of higher order mean curvatures $H_k$. We show that under a suitable pinching condition, the hypersurface is diffeomorpic and almost isometric to a standard sphere. Moreover, as a corollary, we show that a hypersurface of the Euclidean space which is almost Einstein is diffeomorpic and almost isometric to a standard sphere.
\end{abstract}
\section{Introduction and Preliminaries}

Let $(M^n,g)$ be a $n$-dimensional compact, connected, oriented manifold without boundary, isometrically immersed by $\phi$ into the $(n+1)$-dimensional Euclidean space, $(\R^{n+1},can)$, {\it i.e.}, $\phi^*can=g$. If, in addition, $(M^n,g)$ is Einstein, then a well-known result of Cartan and Thomas (\cite{Tho}), also proved by Fialkow (\cite{Fia}), says that $M$ is a round sphere $\Ss^n(R)$ of corresponding radius. A natural question is to ask what one could say if $(M^n,g)$ is almost-Einstein, that is, $||\Ric-kg||_{\infty}\leqslant\varepsilon$, for some positive constant $k$.\\
\indent
Recently, J.F. Grosjean gave a new proof of the Thomas-Cartan theorem using an upper bound of the first eigenvalue of the Laplacian. Indeed, Grosjean proved in \cite{Gr1} that if $(M^n,g)$ has positive scalar curvature, then the first eigenvalue of the Laplacian satisfies
$$\lambda_1(M)\leqslant\frac{1}{n-1}||\sscal||_{\infty},$$
with equality only for geodesic spheres (here $\sscal$ denotes the scalar curvature).\\
\indent
If $(M^n,g)$ is Einstein, {\it i.e.}, $\Ric=(n-1)g$, we know by the Lichnerowicz theorem that $\lambda_1(M)\geqslant n$, and by the above upper bound
$$\lambda_1(M)\leqslant\frac{1}{n-1}||\sscal||_{\infty}=n.$$
So $\lambda_1(M)=n$ and we are in the equality case of both inequalities, that is, $M=\Ss^n$.\\
\indent
This approach leads naturally to consider a pinching result on the first eigenvalue of the Laplacian, which allows to show that an almost Einstein hypersurface of $\R^{n+1}$ is close to a sphere.\\
\indent
First, we can deduce from a theorem of Aubry (\cite{Aub}), which is a pinching theorem corresponding to the Lichnerowicz inequality, that if $\varepsilon$ is small enough, then $M$ is homeomorphic to $\Ss^n$ (see Theorem \ref{corollaireaubry}).\\
\indent
Nethertheless, Aubry's result does not yield to a sufficiently strong rigidity result. For this, we will study another pinching of the first eigenvalue of the Laplacian, which is associated with an extrinsic upper bound involving the scalar curvature. In fact, in this paper, we are interested in more general upper bounds in terms of higher order mean curvatures.\\
\indent 
In \cite{Re2}, Reilly gives upper bounds for $\lambda_1(M)$, in terms of higher order mean curvatures $H_k$, which are defined as the symmetric polynomials in the principal curvatures. He shows that for all $k\in\{1,\cdots,n\}$:
\beqt\label{reilly2}
\lambda_1(M)\left( \int_MH_{k-1}\right)^2 \leqslant\frac{n}{V(M)}\int_MH_k^2,
\eeqt 
with equality only for the standard spheres of $\R^{n+1}$.\\
\indent
By the H\"older inequality, we get similar inequalities with the $L^{2p}$-norms \\$(p\geqslant1)$ of $H_k$:
\beqt\label{reilly3}
\lambda_1(M)\left( \int_MH_{k-1}\right)^2 \leqslant\frac{n}{V(M)^{1/p}}||H_k||_{2p}^2.
\eeqt 
As for inequality (\ref{reilly2}), the equality case in (\ref{reilly3}) characterizes the standard spheres.\\
\indent
Then, a natural question is to know if there exists a pinching result as the following theorem proved by B. Colbois and J.F. Grosjean (\cite{CG})? {\it For $p\geqslant2$ and any $\varepsilon>0$, there exists a constant $C_{\varepsilon}$ depending only on $n$ and $\Hinf$ so that if the pinching condition
$$\frac{n}{V(M)^{1/p}}\HHp-C_{\varepsilon}<\lambda_1(M)$$
is true, then the Haussdorff distance between $M$ and the sphere $S\left( 0,\sqrt{\frac{n}{\lambda_1(M)}}\right) $ is at most $\varepsilon$.}\\
\indent
We give a positive answer to this question, and, as we will see, the case $k=2$, that is involving $H_1$ and $H_2$, will solve the problem for almost-Einstein hypersurfaces.
\begin{thm}\label{thm1}
 Let $(M^n,g)$ be a compact, connected, oriented Riemannian manifold without boundary isometrically immersed in $\R^{n+1}$ and $p_0$ the center of mass of $M$. Assume that $V(M)=1$ and let $k\in\{1,\cdots,n\}$ such that $H_k>0$. Then, for any $p\geqslant2$ and for any $\varepsilon>0$, there exists a constant $C_{\varepsilon}$ depending only on $\varepsilon$, $n$, $\Hinf$ and $\Hkp$ such that
\beqt\tag{$P_{C_{\varepsilon}}$}
\lambda_1(M)\left( \int_MH_{k-1}\right)^2 -\frac{n}{V(M)^{1/p}}||H_k||_{2p}^2>-C_{\varepsilon}
\eeqt
is satisfied, then
\be[i)]
\item $\phi(M)\subset B\left( x_0,\sqrt{\frac{n}{\lambda_1(M)}}+\varepsilon\right)\setminus B\left( x_0,\sqrt{\frac{n}{\lambda_1(M)}}-\varepsilon\right)$.
\item $\forall x\in S\left( x_0,\sqrt{\frac{n}{\lambda_1(M)}}\right),\ B(x,\varepsilon)\cap\phi(M)\neq\emptyset$.
\ee
\end{thm}
We recall that the Haussdorff distance between two compact subsets $A$ and $B$ of a metric space $(E,d)$ is given by
$$d_H(A,B)=\iinf\left\lbrace \eta>0\big|B\subset V_{\eta}(A)\ \text{and}\ a\subset V_{\eta}(B)\right\rbrace,$$
where for any subset $A$, the set $V_{\eta}(A)$ is the tubular neighborhood of $A$ defined by $V_{\eta}(A)=\left\lbrace x\in E\big|d(x,A)<\eta\right\rbrace$. So, $i)$ and $ii)$ of Theorem \ref{thm1} imply that the Haussdorff distance between $M$ and $S\left( x_0,\sqrt{\frac{n}{\lambda_1(M)}}\right)$ is at most $\varepsilon$.
\begin{rem}
We will see in the proof that $C_{\varepsilon}\lgra0$ when $\Hinf\lgra\infty$ or $\varepsilon\lgra0$.
\end{rem}
In this second theorem, if the pinching condition is strong enough, with a control on the $L^{\infty}$-norm of the second fundamental form $B$, we obtain that $M$ is diffeomorphic and almost-isometric to a round sphere in the following sense
\begin{thm}\label{thm2}
 Let $(M^n,g)$ be a compact, connected, oriented Riemannian manifold without boundary isometrically immersed in $\R^{n+1}$ and $p_0$ the center of mass of $M$.  Assume that $V(M)=1$ and let $k\in\{1,\cdots,n\}$ such that $H_k>0$. Then for any $p\geqslant2$, there exists a constant $K$ depending only on $n$, $\Binf$ and $\Hkp$ such that if the pinching condition
 \beqt\tag{$P_K$}
\lambda_1(M)\left( \int_MH_{k-1}\right)^2 -\frac{n}{V(M)^{1/p}}||H_k||_{2p}^2>-K
\eeqt
 is satisfied, then $M$ is diffeomorphic to $\Ss^n$.\\
\indent
More precisely, there exists a diffeomorphism $F$ from $M$ into the sphere $\Ss^n\left( \sqrt{\frac{n}{\lambda_1(M)}}\right) $ of radius  $\sqrt{\frac{n}{\lambda_1(M)}}$ which is a quasi-isometry. Namely, for any $\theta\in]0,1[$, there exists a constant $K_{\theta}$ depending only on $\theta$, $n$, $\Binf$ and $\Hkp$ so that the pinching condition with $K_{\theta}$ implies
$$\big||dF_x(u)|^2-1\big|\leq\theta,$$
for any unitary vector $u\in T_xM$.
\end{thm}
\begin{rem}
We will see in the proof that the constants $C_{\varepsilon}$, $K$ and $K_{\theta}$ of Theorems \ref{thm1} and \ref{thm2} do not depend on $\Hkp$ if $p\geqslant\frac{n}{2k}$.
\end{rem}
These results have a double interest. First, they improve the results in \cite{CG}. Second, the case $k=2$ is especially interesting. Indeed, for hypersurfaces of the Euclidean space, the second mean curvature $H_2$ is, up to a multiplicative constant, the scalar curvature. Precisely, $H_2=\frac{1}{n(n-1)\sscal}$. Then we deduce from Theorems \ref{thm1} and \ref{thm2} two corollaries for almost-Einstein hypersurfaces.\\\\
\indent
Now, we give some preliminaries for the proof of these theorems. Throughout this paper, we consider a compact, connected, oriented Riemannian manifold isometrically immersed in $(\R^{n+1},can)$ by $\phi$. Let $\nu$ be the outward normal unitary vector field. We denote respectively by $\nabla$ and $\Delta$ the Riemannian connection and the Laplacian of $M$, and by $\nablab$ the Riemannian connection of $\R^{n+1}$. Finally, we denote by $\left\langle \cdot,\cdot\right\rangle$ the Euclidean scalar product of $\R^{n+1}$.\\
\indent
The second fundamental $B$ of the immersion is defined by $$B(Y,Z)=\left\langle \nablab_Y\nu,Z\right\rangle$$ and the mean curvature $H$ by $$H=\frac{1}{n}\trace(B).$$
\indent
Now we recall the following well-known identity
\beqt\label{DeltaX}
\frac{1}{2}\Delta|X|^2=nH\left\langle \nu,X\right\rangle-n,
\eeqt
where $X$ is the position vector.\\
\indent
We finish by recalling the definition of higher order mean curvatures. They are extrinsic geometric invariants defined from the second fundamental form and generalizing the mean curvature. We saw that   
$$H=\frac{1}{n}\sigma_1(\kappa_1,\cdots,\kappa_n),$$
where $\sigma_1$ is the first symmetric polynomial and $\kappa_1,\cdots,\kappa_n$ the principal curvatures of $M$. The higher order mean curvatures are defined for $k\in\{1,\cdots,n\}$ by
$$H_k=\frac{1}{\left( \begin{array}{c}
 n\\k 
\end{array}\right) }\sigma_k(\kappa_1,\cdots,\kappa_n),$$
where $\sigma_k$ is the $k$-th symmetric polynomial, that is,
$$\sigma_k(x_1,\cdots,x_n)=\sum_{1\leqslant i_1,\cdots,i_k\leqslant n}x_{i_1}\cdots x_{i_k}.$$
This definition is equivalent to  
$$H_k=\frac{1}{k!\left( \begin{array}{c}n\\k \end{array}\right)}\sum_{\begin{array}{c}1\leqslant i_1,\cdots,i_k\leqslant n\\
1\leqslant j_1,\cdots,j_k\leqslant n\end{array}}\epsilon\left( \begin{array}{c} i_1,\cdots,i_k\\
j_1,\cdots,j_k\end{array}\right) B_{i_1j_1}\cdots B_{i_kj_k},$$
where the $B_{ij}$'s are the coefficients of the second fundamental form $B$ in a local orthonormal frame $\base$. Moreover, we denote by $\epsilon\left( \begin{array}{c} i_1,\cdots,i_k\\
j_1,\cdots,j_k\end{array}\right)$ the usual symbols for permutation. Finally, by convention, we set $H_0=1$ and $H_{n+1}=0$.\\
\indent
These mean curvatures satisfy some properties as the Hsiung-Minkowski formula (see \cite{Hsi})
\beqt\label{Hsuing}
\int_M\big( H_{k-1}-H_k\left\langle X,\nu\right\rangle \big)=0,
\eeqt
and the following inequalities
\begin{lem}\label{inegalitesentreHk} 
If $k\in\{1,\cdots,n\}$, and $H_k$ is a positive function, then
$$H_k^{\frac{1}{k}}\leq H_{k-1}^{\frac{1}{k-1}}\leq\cdots\leq
H_2^{\frac{1}{2}}\leq H.$$
\end{lem}
\section{An $L^2$-approach to the problem}

We prove Theorems \ref{thm1} and \ref{thm2} in two steps. First, we prove that if the pinching condition $(P_C)$ is satisfied, then $M$ is close to a sphere in an $L^2$-sense. For this, we prove a first lemma which gives an estimate of the $L^2$-norm of the position vector $X$.
\begin{lem}\label{lem1}
If the pinching condition $(P_C)$ is satisfied for $C<\frac{n}{2}\HHkp$, then
$$\frac{n\lambda_1(M)\left( \int_MH_{k-1}\right)^4 }{\left( C+\lambda_1(M)\left( \int_MH_{k-1}\right)^2\right)^2 }\leqslant||X||_2^2\leqslant\frac{n}{\lambda_1(M)}\leqslant A_1,$$
where $A_1$ is a positive constant depending only on $n$, $\Hinf$ and $\Hkp$.
\end{lem}
\noindent
{\it Proof:} If $(P_C)$ is satisfied, we have:
$$
\lambda_1(M)\left( \int_MH_{k-1}\right)^2 \geqslant n\HHkp-C.
$$
If, in addition, we assume that $C<\frac{n}{2}\Hkp$, we get
$$\lambda_1(M)\left( \int_MH_{k-1}\right)^2\geqslant\frac{n}{2}\HHkp,$$
and so
\beqt\label{encadrement1}
\frac{n}{\lambda_1(M)}\leqslant\frac{2\left( \int_MH_{k-1}\right)^2}{\HHkp}\leqslant\frac{2||H||_{\infty}^{2(k-1)}}{\HHkp}.
\eeqt
Moreover, by the variational characterization of $\lambda_1(M)$, we have
$$\lambda_1(M)\int_M|X|^2\leqslant\int_M\left( \sum_{i=1}^{n+1}|dX_i|^2\right)=n,$$
where $X_i$ are the functions defined by $X=\sum_{i=1}^{n+1}X_i\partial_i$, where $\{\partial_1,\cdots,\partial_{n+1}\}$ is the canonical frame of $\R^{n+1}$. So we have $||X||_2^2\leqslant\dfrac{n}{\lambda_1(M)}$, and by (\ref{encadrement1}),
$$||X||_2^2\leqslant A_1(n,\Hinf,\Hkp).$$
For the left hand side, we have
\beQ
\lambda_1(M)\left( \int_M|X|^2\right) \left( \int_MH_{k-1}\right)^4&\leqslant&n \left( \int_MH_{k-1}\right)^4\\
&\leqslant&n\left( \int_MH_{k}\left\langle X,\nu\right\rangle \right)^4\\
&\leqslant&n\left( \int_MH_k^2\right)^2 \left( \int_M|X|^2\right)^2 .
\eeQ
Then, by the H\"older inequality, we deduce
$$\lambda_1(M)\left( \int_MH_{k-1}\right)^4\leqslant n\HHkp \left( \int_M|X|^2\right),$$
and with the pinching condition, 
$$||X||_2^2\geqslant\frac{n\lambda_1(M)\left( \int_MH_{k-1}\right)^4 }{\left( C+\lambda_1(M)\left( \int_MH_{k-1}\right)^2\right)^2 }.$$
$\finpreuve$
From now, we denote by $X^T$ the orthogonal projection of $X$ on $M$. That is, if for $x\in M$, $\base$ is an orthonormal frame of $T_xM$, then 
$$X^T=\sum_{i=1}^n\left\langle  X,e_i\right\rangle e_i=X-\left\langle  X,\nu\right\rangle\nu.$$ 
In the following lemma, we show that the pinching condition $(P_C)$ implies that the $L^2$-norm of $X^T$ is close to zero.
\begin{lem}\label{lem2}
The pinching condition $(P_C)$ with $C<\frac{n}{2}\HHkp$ implies 
$$||X^T||_2^2\leqslant A_2C,$$
where $A_2$ is a positive constant depending only on $n$, $\Hinf$ and $\Hkp$.
\end{lem}
\noindent{\it Proof:} We saw that
\beQ
\lambda_1(M)\int_M|X|^2&\leqslant&n,
\eeQ
so by the Hsiung-Minkowski formula and the Cauchy-Schwarz inequality
\beQ
\lambda_1(M)\int_M|X|^2\left( \int_MH_{k-1}\right)^2 &\leqslant&n\left( \int_MH_{k-1}\right)^2\\
&\leqslant& n\left( \int_MH_k\left\langle X,\nu\right\rangle \right)^2\\
&\leqslant&||H_k||_2^2\int_M\left\langle  X,\nu\right\rangle^2\\
&\leqslant&||H_k||_{2p}^2\int_M\left\langle  X,\nu\right\rangle^2.\\
\eeQ
Then we deduce 
\beQ
n||H_k||_{2p}^2||X^T||_2^2&\leqslant&n\HHkp\left( \int_M\left( |X|^2-\left\langle X,\nu\right\rangle^2 \right) \right) \\
&\leqslant&n\HHkp\left[ \int_M|X|^2-\lambda_1(M)\left( \int_MH_{k-1}\right)^2\int_M|X|^2\right]\\
&\leqslant&\left[ n\HHkp-\lambda_1(M)\left( \int_MH_{k-1}\right)^2\right]||X||_2^2\\
&\leqslant&C||X||_2^2\leqslant A_1C.
\eeQ
Finally, we get
$$||X^T||_2^2\leqslant\frac{A_1C}{n||H_k||_{2p}^2}=A_2C.$$
$\finpreuve$\\
\indent
In order to prove assertion $i)$ of Theorem \ref{thm1}, we will show that 
$$\bigg|\bigg||X|-\sqrt{\frac{n}{\lambda_1(M)}} \bigg|\bigg|_{\
\infty}\leqslant\varepsilon.$$
 For this, we need an upper bound on the $L^2$-norm of the function 
 $$\varphi:=|X|\left( |X|-\sqrt{\frac{n}{\lambda_1(M)}}\right)^2.$$ 
Before getting such an estimate, we introduce the two following vector fieds:
$$\left\lbrace 
\begin{array}{l}
Y=nH_k\nu-\lambda_1(M)\left( \int_MH_{k-1}\right)X,  \\ \\
Z=\sqrt{\dfrac{n}{\lambda_1}}\dfrac{|X|^{1/2}H_k}{\left( \int_MH_{k-1}\right)} \nu-\dfrac{X}{|X|^{1/2}}.
\end{array}\right. $$
First, we have the following:
\begin{lem}\label{lem3}
The pinching condition $(P_C)$ implies 
$$||Y||_2^2\leqslant nC.$$
\end{lem}
\noindent
{\it Proof:} We have
\beQ
||Y||_2^2&=&n^2\int_MH_k^2+\lambda_1(M)\left( \int_MH_{k-1}\right)^2\int_M|X|^2\\
&&-2n\lambda_1(M)\left( \int_MH_{k-1}\right)\int_MH_k\left\langle X,\nu\right\rangle\\
&\leqslant&n^2\HHkp+n\lambda_1(M) \left( \int_MH_{k-1}\right)^2-2n\lambda_1(M)\left( \int_MH_{k-1}\right)^2\\
&\leqslant&n\left( n\HHkp-n\lambda_1(M) \left( \int_MH_{k-1}\right)^2\right) \\
&\leqslant&nC, 
\eeQ
where we used the Hsiung-Minkowski formula (\ref{Hsuing}), and the fact that 
$$||X||_2^2\leqslant\frac{n}{\lambda_1(M)}.$$
$\finpreuve$
We also have
\begin{lem}\label{lem4}
If the pinching condition $(P_C)$ is satisfied, with $C<\frac{n}{2}\HHkp$, then
$$||Z||_2^2\leqslant A_3C,$$
where $A_3$ is a positive constant depending only on $n$, $\Hinf$ and $\Hkp$.
\end{lem}
\noindent
{\it Proof:} We have
\beQ
||Z||_2^2&=&\frac{n}{\lambda_1(M)\left( \int_MH_{k-1}\right)^2 }\int_M|X|H_k^2+\int_M|X|-2\frac{\sqrt{\frac{n}{\lambda_1(M)}}}{\int_MH_{k-1}}\int_MH_k\left\langle X,\nu\right\rangle \\
&\leqslant&\frac{n}{\lambda_1(M)\left( \int_MH_{k-1}\right)^2 }\int_M|X|H_k^2+\int_M|X|-2\sqrt{\frac{n}{\lambda_1(M)}}
\eeQ
By the H\"older inequality, we get
\beQ
||Z||_2^2&\leqslant&\frac{n}{\lambda_1(M)\left( \int_MH_{k-1}\right)^2 }\left( \int_MH_k^4\right)^{1/2} \left( \int_M|X|^2\right)^{1/2}\\
&&+\left( \int_M|X|^2\right)^{1/2} -2\sqrt{\frac{n}{\lambda_1(M)}}\\
&\leqslant&\sqrt{\frac{n}{\lambda_1(M)}}\left[ \frac{n}{\lambda_1(M)}\frac{\HHkp}{\left( \int_MH_{k-1}\right)^2} -1\right] \\
&\leqslant&\left( \frac{n}{\lambda_1(M)}\right) ^{3/2}\frac{1}{\left( \int_MH_{k-1}\right)^2}\left[ n\HHkp-\lambda_1(M)\left( \int_MH_{k-1}\right)^2\right] \\
&\leqslant&A_3C,
\eeQ
where $A_3$ depends only on $n$, $\Hinf$ and $\Hkp$. Note that we have used Lemma \ref{inegalitesentreHk} and the fact that $\dfrac{n}{\lambda_1(M)}\leqslant\dfrac{2||H||_{\infty}^{2(k-1)}}{\HHkp}.$ $\finpreuve$\\
\indent
Now, we give an upper bound for the $L^2$-norm of the function $\varphi$.
\begin{lem}\label{lem5}
The pinching condition $(P_C)$ with $ C<\frac{n}{2}\HHkp$ implies
$$||\varphi||_2\leqslant A_4||\varphi||_{\infty}^{3/4}C^{1/4}.$$
\end{lem}
\noindent
{\it Proof:} We have
$$
||\varphi||_2=\left( \int_M\varphi^{3/2}\varphi^{1/2}\right)^{1/2}\leqslant||\varphi||_{\infty}^{3/4}\big|\big|\varphi^{1/2}\big|\big|_1^{1/2}. 
$$
Moreover,
\beQ
\int_M\varphi^{1/2}&=&\Bigg|\Bigg||X|^{1/2}X-\sqrt{\frac{n}{\lambda_1(M)}}\frac{X}{|X|^{1/2}}\Bigg|\Bigg|_1\\
&=&\Bigg|\Bigg|-\frac{|X|^{1/2}}{\lambda_1(M)\int_MH_{k-1}}Y+\frac{n}{\lambda_1(M)}\frac{|X|^{1/2}H_k}{\int_MH_{k-1}}\nu-\sqrt{\frac{n}{\lambda_1(M)}}\frac{X}{|X|^{1/2}}\Bigg|\Bigg|_1\\
&\leqslant&\Bigg|\Bigg|\frac{|X|^{1/2}}{\lambda_1(M)\int_MH_{k-1}}Y\Bigg|\Bigg|_1+\sqrt{\frac{n}{\lambda_1(M)}}||Z||_1.
\eeQ
By the H\"older inequality, we get
\beQ
\Bigg|\Bigg|\frac{|X|^{1/2}}{\lambda_1(M)\int_MH_{k-1}}Y\Bigg|\Bigg|_1&\leqslant&\frac{1}{\lambda_1}\left( \int_M|X|^2\right)^{1/4}||Y||_2\\
&\leqslant&\frac{A_1^{3/4}}{n^{1/2}}C^{1/2}.
\eeQ
Finally, from Lemmas \ref{lem3} and \ref{lem4}, we obtain
$$||\varphi^{1/2}||_1^{1/2}\leqslant A_4C^{1/4},$$
where $A_4$ is a positive constant depending only on $n$, $\Hinf$ and $\Hkp$. \\
$\finpreuve$
\section{Proof of Theorem \ref{thm1}}

The proof of Theorem \ref{thm1} is an immediate consequence of the two following lemmas:
\begin{lem}\label{lem6}
For $p\geqslant2$ and any $\eta>0$, there exists $K_{\eta}(n,\Hinf,\Hkp)$ so that if $(P_{K_{\eta}})$ is true, then $||\varphi||_{\infty}\leqslant\eta$. Moreover, $K_{\eta}\lgra0$ when $\Hinf\lgra\infty$ or $\eta\lgra0$.
\end{lem}
\begin{lem}\label{lemgeom1}[Colbois-Grosjean \cite{CG}]
Let $x_0$ be a point of the sphere $S(0,R)$ in $\R^{n+1}$ with the center at the origin and of radius $R$. Assume that $x_0=Re$ with $e\in\Ss^n$. Now let $(M^n,g)$ be a compact, connected, oriented n-diemsnional Riemannian manifold without boundary isometrically mmersed in $\R^{n+1}$. If the image of $M$ is contained in $\Big( B(0,R+\eta)\setminus B(0,R-\eta)\Big)\setminus B(x_0,\rho)$ with $\rho=4(2n-1)\eta$. Then there exists a point $y_0\in M$ so that the mean curvature of $M$ in $y_0$ satisfies $|H(y_0)|\geqslant\frac{1}{4n\eta}$.
\end{lem}
We will prove Lemma \ref{lem6} in Section \ref{lemtech}. Now, we will prove Theorem \ref{thm1} by using these two lemmas.
\paragraph{Proof of Theorem \ref{thm1}}
Let $\varepsilon>0$ and consider the function 
$$f(t):=t\left( t-\sqrt{\frac{n}{\lambda_1(M)}}\right)^2 .$$
 We set 
$$\eta(\varepsilon):=\iinf\left\lbrace f\left( \sqrt{\frac{n}{\lambda_1(M)}}-\varepsilon\right) ,f\left( \sqrt{\frac{n}{\lambda_1(M)}}+\varepsilon\right),\frac{1}{27\Hinf^3}\right\rbrace .$$
By definition, $\eta(\varepsilon)>0$, and by Lemma \ref{lem6}, there exists $K_{\eta(\varepsilon)}$ such that for all $x\in M$, 
\beqt\label{f}
f(|X|(x))\leqslant\eta(\varepsilon).
\eeqt
Now to prove the theorem, it is sufficient to assume $\varepsilon<\frac{2}{3\Hinf}$. We will show that either 
\beqt\label{alternative}
\sqrt{\frac{n}{\lambda_1(M)}}-\varepsilon\leqslant|X|\leqslant\sqrt{\frac{n}{\lambda_1(M)}}+\varepsilon\quad\text{or}\quad |X|<\frac{1}{3}\sqrt{\frac{n}{\lambda_1(M)}}
\eeqt
By examining the function $f$, it is easy to see that $f$ has a unique local maximum at $\frac{1}{3}\sqrt{\frac{n}{\lambda_1(M)}}$. Moreover, from the definition of $\eta(\varepsilon)$, we have 
$$\eta(\varepsilon)<\frac{4}{27\Hinf^3}\leqslant\frac{4}{27}\left( \frac{n}{\lambda_1(M)}\right)^{3/2}=f\left( \frac{1}{3}\sqrt{\frac{n}{\lambda_1(M)}}\right) .$$
Since we assume $\varepsilon<\frac{2}{3\Hinf}\leqslant\frac{2}{3}\sqrt{\frac{n}{\lambda_1(M)}}$, we have
$$\frac{1}{3}\sqrt{\frac{n}{\lambda_1(M)}}<\sqrt{\frac{n}{\lambda_1(M)}}-\varepsilon,$$ which with (\ref{f}) yields (\ref{alternative}).\\\\
\indent
Now, from Lemma \ref{lem1}, we deduce that there exists a point $y_0\in M$ such that
$$|X(y_0)|^2\geqslant\frac{n\lambda_1(M)\left( \int_MH_{k-1}\right)^4 }{\left( K_{\eta(\varepsilon)}+\lambda_1(M)\left( \int_MH_{k-1}\right)^2\right)^2 }.$$
Since $K_{\eta(\varepsilon)}<\frac{n}{2}\HHkp$, the condition $(P_C)$ implies 
$$K_{\eta(\varepsilon)}<\frac{n}{2}\HHkp\leqslant\lambda_1(M)\left( \int_MH_{k-1}\right)^2\leqslant2\lambda_1(M)\left( \int_MH_{k-1}\right)^2.$$
We deduce that $$|X(y_0)|\geqslant \frac{1}{3}\sqrt{\frac{n}{\lambda_1(M)}}.$$
Since $M$ is connected, for any $x\in M$,
$$\sqrt{\frac{n}{\lambda_1(M)}}-\varepsilon\leqslant|X|(x)\leqslant\sqrt{\frac{n}{\lambda_1(M)}}+\varepsilon,$$
which proves the assertion $i)$ of the theorem.\\
\indent
In order to prove the second, we consider the pinching condition $(P_{C_{\varepsilon}})$ with $C_{\varepsilon}=K_{\eta\left(\frac{\varepsilon}{4(2n-1)} \right)}$. Then assertion $i)$ is still valid.\\
Let $x= \sqrt{\frac{n}{\lambda_1(M)}}\,e\in S\left( 0,\sqrt{\frac{n}{\lambda_1(M)}}\right) $, with $e\in\Ss^n$ and assume that\\ $B(x,\varepsilon)\cap M=\emptyset.$ We can apply Lemma \ref{lemgeom1}. So, there exists a point $y_0\in M$ such that $|H(y_0)|\geqslant\frac{2n-1}{n\varepsilon} >\Hinf$ since we assumed $\varepsilon<\frac{2}{3\Hinf}\leqslant\frac{2n-1}{n\Hinf}$. This is a contradiction and so $B(x,\varepsilon)\cap M\neq\emptyset$. The assertion $ii)$ is satisfied and $C_{\varepsilon}\lgra0$ when $\Hinf\lgra\infty$ or $\varepsilon\lgra0$. $\finpreuve$
\section{Proof of Theorem \ref{thm2}}

From Theorem \ref{thm1}, we know that for any $\varepsilon>0$, there exists $C_{\varepsilon}$ depending only on $n$, $\Hinf$ and $\Hkp$ so that if $(P_{ C_{\varepsilon}})$ is true, then
$$\bigg||X|(x)-\sqrt{\frac{n}{\lambda_1(M)}}\bigg|\leqslant\varepsilon$$
for all $x\in M$. Since $\sqrt{n}\Hinf\leqslant\Binf$, it is easy to see that we can assume that $C_{\varepsilon}$ depends only on $n$, $\Binf$ and $\Hkp$.\\
\indent
The proof of Theorem \ref{thm2} is an immediate consequence of the following lemma about the $L^{\infty}$-norm of $X^T$.
\begin{lem}\label{lem7}
For $p\geqslant2$ and any $\eta>0$, there exists $K_{\eta}(n,\Binf,\Hkp)$ so that if $(P_{k_{\eta}})$ is true, then $||X^T||_{\infty}\leqslant\eta$.
\end{lem}
We will prove this lemma in Section \ref{lemtech}.
\paragraph{Proof of Theorem \ref{thm2}}
Let $\varepsilon<\frac{1}{2}\sqrt{\frac{n}{\Binf}}\leqslant\sqrt{\frac{n}{\lambda_1(M)}}$. This choice of $\varepsilon$ implies that if the pinching condition $(P_{C_{\varepsilon}})$ is true, then $|X|$ never vanishes, and so we can consider the following map 
\function{F}{M}{S\left( 0,\sqrt{\frac{n}{\lambda_1(M)}}\right)}{x}{\sqrt{\frac{n}{\lambda_1(M)}}\;\frac{X}{|X|}.}
Without any pinching condition, a straightforward computation yields to
\beqt\label{dF}
\Big||dF_x(u)|^2-1\Big|\leqslant \bigg|\frac{n}{\lambda_1(M)}\frac{1}{|X|^2}-1\bigg|+\frac{n}{\lambda_1(M)}\frac{1}{|X|^4}\left\langle u,X\right\rangle^2, 
\eeqt
for any unitary vector $u\in T_xM$. But,
$$ \bigg|\frac{n}{\lambda_1(M)}\frac{1}{|X|^2}-1\bigg|=\frac{1}{|X|^2}\bigg|\frac{n}{\lambda_1(M)}-|X|^2\bigg|\leqslant\varepsilon\frac{\sqrt{\frac{n}{\lambda_1(M)}}+|X|}{|X|^2}\leqslant\varepsilon\frac{2\sqrt{\frac{n}{\lambda_1(M)}}+\varepsilon}{\left( \sqrt{\frac{n}{\lambda_1(M)}}-\varepsilon\right)^2 }$$
We recall that $\frac{n}{A_1}\leqslant\lambda_1\leqslant\Binf^2$. Since we assume $\varepsilon<\frac{1}{2}\sqrt{\frac{n}{\Binf}}$, the right hand side is bounded by a constant depending only on $n$, $\Binf$ and $\Hkp$. So we have
\beqt\label{dF1}
\bigg|\frac{n}{\lambda_1(M)}\frac{1}{|X|^2}-1\bigg|\leqslant\varepsilon\gamma(n,\Binf,\Hkp).
\eeqt
Moreover, since $C_{\varepsilon}\lgra0$ when $\varepsilon\lgra0$, there exists $\varepsilon(n,\Binf,\Hkp,\eta)$ so that $C_{\varepsilon}\leqslant K_{\eta}$ (where $K_{\eta}$ is the constant of Lemma \ref{lem7}) and so, $||X^T||_{\infty}\leqslant\eta$. As before, there exists a constant $\delta$ depending on $n$, $\Binf$ and $\Hkp$ such that 
\beqt\label{dF2}
\frac{n}{\lambda_1(M)}\frac{1}{|X|^4}\left\langle u,X\right\rangle^2\leqslant\frac{n}{\lambda_1(M)}\frac{1}{|X|^4}||X^T||_{\infty}^2\leqslant\eta^2\delta(n,\Binf,\Hkp).
\eeqt
Then, from (\ref{dF}), (\ref{dF1}) and (\ref{dF2}), we deduce that $(P_{C_{\varepsilon}})$ implies
$$\Big||dF_x(u)|^2-1\Big|\leqslant\varepsilon\gamma+\eta^2\delta.$$
Take $\eta=\sqrt{\frac{\theta}{2\delta}}$. We can assume that $\varepsilon$ is small enough to have $\varepsilon\gamma\leqslant\frac{\theta}{2}$. In that case, we have
$$\Big||dF_x(u)|^2-1\Big|\leqslant\theta.$$
Now, it is sufficient to fix $\theta\in]0,1[$ and, $F$ is a local diffeomorphism from $M$ into $S\left( 0,\sqrt{\frac{n}{\lambda_1(M)}}\right)$. Since $S\left( 0,\sqrt{\frac{n}{\lambda_1(M)}}\right)$ is simply connected for $n\geqslant2$, the map $F$ is a global diffeomorphism.$\finpreuve$
\section{Application to almost-Einstein Hypersurfaces}

In this section, we give an application of Theorems \ref{thm1} and \ref{thm2} to almost Einstein hypersurfaces.  In fact, we obtain two different rigidty results.
\begin{cor}\label{cor1}
Let $(M^n,g)$ be a connected, oriented Riemannian manifold without boundary isometrically immersed in $\R^{n+1}$. If $(M^n,g)$ is almost-Einstein, that is, $||\Ric-kg||_{\infty}\leqslant\varepsilon$ for a positive constant $k$, with $\varepsilon$ small enough depending on $n$, $k$ and $\Hinf$, then 
$$d_H\left( M,\Ss^n\left(\sqrt{\frac{n-1}{k}}\right)\right) \leqslant\varepsilon.$$
\end{cor}
\begin{cor}\label{cor2}
Let $(M^n,g)$ be a compact, connected, oriented Riemannian manifold without boundary isometrically immersed in $\R^{n+1}$. If $(M^n,g)$ is almost-Einstein, that is, $||\Ric-kg||_{\infty}\leqslant\varepsilon$ for a positive constant $k$, with $\varepsilon$ small enough depending on $n$, $k$ and $\Binf$, then $M$ is diffeomorphic and almost isometric to $\Ss^n\left(\sqrt{\frac{n-1}{k}}\right)$
\end{cor}\noindent
{\it Proof:} Assume that $k=n-1$. By the assumption $||\Ric-(n-1)g||_{\infty}\leqslant\varepsilon$, the Lichnerowicz theorem implies that
$$\lambda_1(M)\geqslant \frac{n(n-1-\varepsilon)}{n-1}=n-\frac{n\varepsilon}{n-1}.$$
So, for $p\geqslant2$, we have
\beQ
\lambda_1(M)\left( \int_MH\right)^2-n||H_2||_{2p}^2 &\geqslant& n\left( 1-\frac{\varepsilon}{n-1}\right) \left( \int_MH_2^{1/2}\right)^2-n ||H_2||_{2p}^2\\
&\geqslant& n\left( 1-\frac{\varepsilon}{n-1}\right) \iinf\{H_2\}-n\ssup\{H_2\}\\
&\geqslant& n\left( 1-\frac{\varepsilon}{n-1}\right)^2-n \left( 1+\frac{\varepsilon}{n-1}\right)\\
&\geqslant&\frac{n\varepsilon^2}{(n-1)^2}-\frac{3n\varepsilon}{n-1}=-\beta_n(\varepsilon),
\eeQ
where $\beta_n$ is a positive function such that $\beta_n(\varepsilon)\lgra0$ when $\varepsilon\lgra0$.\\
\indent
We can choose $\varepsilon$ small enough to have $\beta_n(\varepsilon)\leqslant C(n,\Binf,||H_2||_{2p})$ of Theorem \ref{thm2}, and we deduce that there exists $\varepsilon$ depending only on $n$ and $\Binf$ so that if $||\Ric-(n-1)||_{\infty}\leqslant\varepsilon$, then $M$ is diffeomorphic and almost isometric to $\Ss^n$. Since $1-\frac{n}{\varepsilon}\leqslant H_2\leqslant1+\frac{n}{\varepsilon}$, there is no dependence on $||H_2||_{2p}$. By homothety, we get the result for any $k>0$.\\
\indent
The proof or Corollary \ref{cor1} is the same, we use Theorem \ref{thm1} instead of Theorem \ref{thm2}.
$\finpreuve$\\
\indent
As we mentioned, these two corollaries are to be compared to the following thoerem obtained by a pinching result associated with the Lichnerowicz inequality.
\begin{thm}\label{corollaireaubry}
Let $(M^n,g)$ be a compact, connected, oriented Riemannian manifold without boundary isometrically immersed in $\R^{n+1}$ and $p>\frac{n}{2}$. Then for any $k>0$, there exists $\varepsilon(k,n,||K||_{2p})$ (where $K$ is the sectional curvature of $M$) such that if 
$$||\Ric-kg||_{\infty}\leqslant\varepsilon,$$
then $M$ is homeomorphic to $\Ss^n$.
\end{thm}
\noindent
{\it Proof:} The assumption 
$$||\Ric-kg||_{\infty}\leqslant\varepsilon,$$
implies that the scalar curvature satisfies
$$0<n(k-\varepsilon)\leqslant\sscal\leqslant n(k+\varepsilon).$$
So the first eigenvalue of the Laplacian can be bounded form above
$$\lambda_1(M)\leqslant\frac{n(k+\varepsilon)}{n-1}.$$
On the other hand, the Lichnerowicz theorem says that
$$\lambda_1(M)\geqslant \frac{n(k-\varepsilon)}{n-1}.$$
Now, let us recall the following theorem due to E. Aubry (\cite{Aub}), which a generalization of a theorem of Ilias (\cite{Ili}).
\begin{thrm}[Aubry \cite{Aub}]
Let $p$, $R$ and $A$ be some real numbers such that $p>\frac{n}{2}$, $R>0$ and $A>0$. Let $(M^n,g)$ be a complete Riemannian manifold. There exists $\alpha(p,n,A)>0$ such that if $$\ssup_x\frac{\big|\big|\big(\underline{\Ric}-(n-1)\big)^-\big|\big|_{L^p(B(x,R))}}{V(B(x,R))}\leqslant\alpha(p,n,A),$$ $||K||_{2p}\leqslant A$, and $$\lambda_1(M)\leqslant n(1+\alpha(p,n,A)),$$ then $M$ is homeomorphic to $\Ss^n$.
\end{thrm}
In this theorem, $\underline{\Ric}(x)$ is the smallest eigenvalue of the symmetric bilinear form $\Ric(x)$ on $T_xM$, and $ \big(\underline{\Ric}-(n-1)\big)^-={\rm{max}}\big(0,-\underline{\Ric}+(n-1)\big)$.\\
\indent
Since $M$ is almost-Einstein, we are precisely in the assumptions of this theorem, and it is sufficient to choose $\varepsilon(k,n,||K||_{2p})>0$ small enough.
$\finpreuve$

\section{Proof of the technical Lemmas}
\label{lemtech}
The proof of Lemmas \ref{lem6} and \ref{lem7} is based on the following result due to Colbois and Grosjean \cite{CG} using a Niremberg-Moser type of argument.
\begin{lem}\label{lemfond}
Let $(M^n,g)$ be a compact, connected, oriented Riemannian manifold without boundary isometrically immersed in $\R^{n+1}$ and let $\xi$ be a nonnegative continuous function on $M$ such that $\xi^k$ is smooth for $k\geqslant2$. Assume there exist $0\leqslant l<m\leqslant2$ such that
$$\frac{1}{2}\xi^{2k-2}\Delta\xi^2\leqslant\ddiv\omega+(\alpha_1+k\alpha_2)\xi^{2k-l}+(\beta_1+k\beta_2)\xi^{2k-m},$$
where $\omega$ is a $1$-form and $\alpha_1$, $\alpha_2$, $\beta_1$ and $\beta_2$ are some nonnegative constants. Then for any $\eta>0$, there exists a positive constant $L$ depending on $n$, $\delta$, $\alpha_1$, $\alpha_2$, $\beta_1$, $\beta_2$, $\Hinf$ and $\eta$ so that if $||\xi||_{\infty}>\eta$, then
$$||\eta||_{\infty}\leqslant L||\xi||_2.$$
Moreover, $L$ is bounded when $\eta\lgra+\infty$ and if $\beta_1>0$, then $L\lgra+\infty$ if $\Hinf\lgra+\infty$ or $\eta\lgra0$.
\end{lem}
In order to prove \ref{lem6} and \ref{lem7}, it is sufficient to find an upper bound for the functions
$$\left\lbrace \begin{array}{l}
 \varphi^{2k-2}\Delta\varphi^2\\ 
|X^T|^{2k-2}\Delta|X^T|^2.
\end{array}\right. $$
For this, the pinching condition $(P_C)$ is used only one time, to obtain an upper bound of $||X||_{\infty}$ depending only on $n$, $\Hinf$ and $\Hkp$. 
\begin{lem}\label{lem8}
If the pinching condition $(P_C)$ is satisfied with $C<\frac{n}{2}\HHkp$, then there exists $E(n,\Hinf,\Hkp)$ such that $||X||_{\infty}\leqslant E$.
\end{lem}
\noindent
{\it Proof:} From (\ref{DeltaX}), we have
$$\frac{1}{2}\Delta|X|^2|X|^{2l-2}\leqslant n\Hinf|X|^{2l-1}.$$
We apply Proposition \ref{lemfond} to the function $\xi=|X|$. We now that if $||X||_{\infty}>E$, then there exists a constant $L(n,\Hinf,E)$ such that 
$$||X||_{\infty}\leqslant L||X||_2.$$
By the pinching condition $(P_C)$ with $C<\frac{n}{2}\HHkp$, we obtain from Lemma \ref{lem1} that 
$$||X||_{\infty}\leqslant LA_1(n,\Hinf,\Hkp)^{1/2}.$$
But, $L$ is bounded when $E\lgra0$, so we can choose $E=E(n,\Hinf,\Hkp)$ big enough to have 
$$LA_1(n,\Hinf,\Hkp)^{1/2}<E.$$
In that case, we have $||X||_{\infty}\leqslant E(n,\Hinf,\Hkp)$.
$\finpreuve$
Then, the proof of Lemmas \ref{lem6} and \ref{lem7} is exactly the same as the proof of the technical Lemmas in \cite{CG}, \cite{Roth} or \cite{Roth2}.
\bibliographystyle{amsplain}
\bibliography{laplacian}
\end{document}